\documentclass[11pt]{amsart}
\usepackage{euler}

\usepackage[hidelinks, colorlinks=true, allcolors=blue]{hyperref}

\input{./commands.sty}

\title{Hypergraph associahedra and compactifications of moduli spaces of points}
\author{Jasper Bown and Javier González Anaya}

\begin{document}

\begin{abstract}
We prove that every Hassett compactification of the moduli space of weighted stable rational curves that admits both a reduction map from the Losev-Manin compactification and a reduction map to projective space is a toric variety, whose corresponding polytope is a hypergraph associahedron (also known as a nestohedron). In addition, we present an analogous result for the moduli space of labeled weighted points in affine space up to translation and scaling. These results are interconnected, and we make their relationship explicit through the concept of ``inflation" of a hypergraph associahedron.
\end{abstract}

\maketitle

\section{Introduction}
Moduli theory is a rich source of interesting combinatorial structures, and it is common for combinatorial insights to shed light on the geometric properties and behavior of moduli spaces. A remarkable case study is offered by moduli spaces of ordered points in varieties; see \cite{ascher2020enumerating,clader2023wonderful,di2022families,gillespie2023degenerations, moon2018birational,kannan2021chow} for recent examples. In this paper, we contribute to this body of work by providing a partial characterization of toric Hassett compactifications of $M_{0,n}$, the moduli space of $n$ distinct labeled points in $\bb{P}^1$. Additionally, we present a parallel result for the moduli space of $n$ distinct labeled points in $d$-dimensional affine space, denoted $\tdno$. These results are interconnected, and we make their relationship explicit through the concept of inflation of hypergraph associahedra.

Hassett's theory \cite{hassett2003moduli} offers a framework for constructing compactifications of $M_{0,n}$ using weight data. Specifically, given a \emph{collection of weight data}, this is, a tuple $\cal{A} = (a_1,\dots,a_n)\in\bb{Q}^n$ with $0 < a_i \leq 1$ and $\sum_{i=1}^n a_i > 2$, Hassett constructs a compactification $\overline{M}_{0,\cal{A}}$ of $M_{0,n}$. Several known compactifications fall into this category, for example the Grothendieck-Knudsen compactification $\overline{M}_{0,n}$ (for $\cal{A}=(1,\dots,1)$), the Losev-Manin space $\overline{M}^{LM}_{0,n}$ (for $\cal{A}_{LM}$ as in Definition~\ref{mon: def LM}), and the compactification given by $\bb{P}^{n-3}$ (for $\cal{A}_{\bb{P}}$ as in Definition~\ref{mon: def pn}).

A natural question then arises: which Hassett compactifications are toric, and which polytopes correspond to these toric Hassett spaces? Ferreira da Rosa, Jensen, and Ranganathan \cite{da2016toric} addressed this by classifying Hassett spaces that are isomorphic to the toric variety of a graph associahedron. Notably, their results include $\overline{M}_{0,n}^{LM}$, which is isomorphic to the toric variety of the permutohedron, itself the graph associahedron of a complete graph. 
Our results extend these findings by employing hypergraph associahedra (also known as nestohedra), a class of polytopes obtained by truncating a simplex based on specific conditions defined by a hypergraph \cite{dovsen2011hypergraph}.

Let $\cal{D}_{0,n}$ denote the set of all collections of weight data. We write $\cal{A}\geq\cal{B}$ if $\cal{B}=(b_1,\dots,b_n)\in\cal{D}_{0,n}$ satisfies $a_i \geq b_i$ for all $i$. When $\cal{A}\geq\cal{B}$, there is a natural reduction map $\overline{M}_{0,\cal{A}} \to \overline{M}_{0,\cal{B}}$. Hassett shows that $\cal{D}_{0,n}$ has a coarse chamber decomposition such that if $\cal{A}\geq\cal{B}$ and $\cal{A}$ and $\cal{B}$ are in the same coarse chamber, then the reduction map is an isomorphism. For convenience let us write $\cal{A}\geq_c\cal{B}$ when there exist collections of weight data $\cal{A}'$ and $\cal{B}'$ in the same coarse chambers as $\cal{A}$ and $\cal{B}$, respectively, such that $\cal{A}'\geq\cal{B}'$.

\begin{theorem}\label{intro: main theorem mon}[Theorem~\ref{mon: main theorem}]
    Let $n\geq 4$ and $\cal{A}\in\cal{D}_{0,n}$ be a collection of weight data such that $\cal{A}_{LM}\geq_c\cal{A}$ and $\cal{A}\geq_c\cal{A}_{\bb{P}}$, so that there are reduction maps $\overline{M}_{0,n}^{LM}\to\overline{M}_{0,\cal{A}}\to\bb{P}^{n-3}$. Then, $\overline{M}_{0,\cal{A}}\cong X(H_\cal{A})$, where $X(H_\cal{A})$ is the toric variety corresponding to the hypergraph associahedron of the hypergraph $H_\cal{A}$.
\end{theorem}

The hypergraph $H_\cal{A}$ is constructed from the weight data $\cal{A}=(a_1,\dots,a_n)$; see Definition~\ref{mon: def HA}. Its vertex set is $[n-2]:=\{1,\dots,n-2\}$, and its hyperedges are the nonempty subsets $I\subseteq[n-2]$ such that $\sum_{i\in I}a_i>1$. The toric variety $X(H_\cal{A})$ is introduced in Definition~\ref{bg: def toric variety of hypergraph}. Besides using known results from \cite{hassett2003moduli,alexeev2008moduli}, our proof also relies on Li's theory of wonderful compactifications \cite{li2009wonderful}.

Our second result concerns a different, yet related, moduli space. This moduli space, denoted $\tdno$, parametrizes collections of $n$ distinct labeled points in $d$-dimensional affine space $\bb{A}_\bb{C}^d$. It was first studied by Chen, Gibney, and Krashen in \cite{Chen-Gibney-Krashen}, where they introduced the compactification $T_{d,n}$. Among their results, they find that $T_{1,n}\cong\overline{M}_{0,n+1}$. In a development analogous to Hassett's compactifications of $M_{0,n}$, Gallardo and Routis \cite{gallardo2017wonderful} constructed many more smooth compactifications by taking collections of weight data and producing modular compactifications.

The set of admissible weights for this moduli problem, denoted by $\cal{D}_{d,n}^T$, consists of all tuples $\cal{A} = (a_1, \dots, a_n)\in\bb{Q}^n$ such that $0 < a_i \leq 1$ and $\sum_{i=1}^n a_i > 1$. The compactification of $\tdno$ corresponding to $\cal{A}\in\cal{D}_{d,n}^T$ is denoted as $T_{d,n}^\cal{A}$. This set can be naturally identified with a wall of $\cal{D}_{0,n+1}$ via the inclusion mapping $\cal{A} = (a_1, \dots, a_n)\in\cal{D}_{d,n}^T$ to $\cal{A}^+ := (a_1, \dots, a_n, 1)\in\cal{D}_{0,n+1}$. Moreover, $\cal{D}_{d,n}^T$ admits a coarse chamber decomposition along which the $T_{d,n}^\cal{A}$ remain constant. This decomposition can be obtained by restricting the fine chamber decomposition of $\cal{D}_{0,n+1}$.

Examples of such compactifications include the Chen-Gibney-Krashen compactification $T_{d,n}$ (for $\cal{A} = (1, \dots, 1)$), the higher-dimensional Losev-Manin compactification $T_{d,n}^{LM}$ (for $\cal{A}_{LM}^T$), a smooth toric compactification introduced in \cite{gallardo2017wonderful} and further studied in \cite{ggg-higher-lm}, and $\bb{P}^{d(n-1)-1}$ (for $\cal{A}_{\bb{P}}^T$). Here  $\cal{A}_{LM}^T,\cal{A}_{\bb{P}}^T\in\cal{D}_{d,n}^T$ are the unique elements such that $(\cal{A}_{LM}^T)^+=\cal{A}_{LM}$ and $(\cal{A}_{\bb{P}}^T)^+=\cal{A}_{\bb{P}}$.

Unlike Hassett's compactifications of $M_{0,n}$, whose explicit construction is not readily available, $T_{d,n}^\cal{A}$ is defined explicitly as an iterated blow-up of $\bb{P}^{d(n-1)-1}$ along a collection of subvarieties $\cal{G}_\cal{A}$. For any $\cal{A}, \cal{B} \in \cal{D}_{d,n}^T$, we write $\cal{A} \geq \cal{B}$ whenever $\cal{A}^+ \geq \cal{B}^+$ , the notation $\cal{A} \geq_c \cal{B}$ is defined in the same way as for $\cal{D}_{0,n}$. It is known that if $\cal{A} \geq_c \cal{B}$, then $\cal{G}_\cal{A} \supseteq \cal{G}_\cal{B}$, and there is a natural reduction morphism $T_{d,n}^{\cal{A}} \to T_{d,n}^{\cal{B}}$ given by the iterated blow-down maps. If $\cal{A}$ and $\cal{B}$ are in the same coarse chamber, then $\cal{G}_\cal{A}=\cal{G}_\cal{B}$ and the reduction map is an isomorphism.

Building on these observations, we note that $T_{d,n}^\cal{A}$ is a smooth toric variety for every $\cal{A}\in\cal{D}_{d,n}^T$ such that $\cal{A}_{LM}^T\geq_c\cal{A}$ and $\cal{A}\geq_c\cal{A}_{\bb{P}}^T$. Our second main result is a characterization of the polytopes defining these toric varieties in terms of our concept of \emph{inflation}. Explicitly, we prove that for any such $\cal{A}$, the polytope defining the toric variety $T_{d,n}^\cal{A}$ is the hypergraph associahedron of the \emph{$d$-inflation} $\operatorname{Inf}_d(H_\cal{A^+})$ of the hypergraph $H_{\cal{A}^+}$ defined above; see Definition~\ref{tdn: def inflation}. Informally, this says that the combinatorial information that goes into these compactifications differs from that one of $\overline{M}_{0,\cal{A}^+}$ simply by ``inflating" the hypergraph. We briefly discuss this concept after the statement of the theorem.

\begin{theorem}\label{intro: main thm tdn}[Theorem~\ref{tdn: main thm tdn}]
    Let $n\geq 4$ and $\cal{A}\in\cal{D}_{d,n}^T$ be a collection of weight data such that $\cal{A}_{LM}^T\geq_c\cal{A}$ and $\cal{A}\geq_c\cal{A}_{\bb{P}}^T$, so that there are reduction maps $T_{0,n}^{LM}\to T_{d,n}^{\cal{A}}\to\bb{P}^{d(n-1)-1}$. Then,
    \[
        T_{d,n}^{\cal{A}} \cong X(\operatorname{Inf}_d(H_{\cal{A}^+})),
    \]
    where $\operatorname{Inf}_d(H_{\cal{A}^+})$ is the $d$-inflation of the hypergraph $H_{\cal{A}^+}$ and $X(\operatorname{Inf}_d(H_{\cal{A}^+}))$ its corresponding toric variety; see Definition~\ref{tdn: def inflation}.
\end{theorem}

\begin{figure}
\centering
\begin{subfigure}{.5\textwidth}
    \centering
    \includegraphics[width=.68\linewidth]{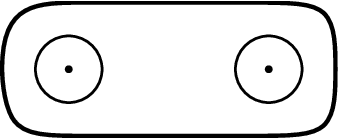}
\end{subfigure}
\begin{subfigure}{.5\textwidth}
    \centering
    \includegraphics[width=.7\linewidth]{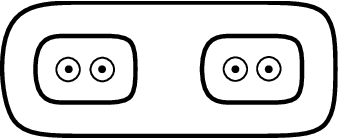}
\end{subfigure}
\caption{The $2$-inflation of a hypergraph with vertex set $\{1,2\}$. The vertex $i$ ``inflates" to the vertices $(i,1)$ and $(i,2)$, and the hyperedge $\{i\}$ inflates to the hyperedge $\{(i,1),(i,2)\}$, for $i=1,2$,.}
\label{fig: intro inflation}
\end{figure}

Here, the vertex set of the $d$-inflation $\operatorname{Inf}_d(H)$ of a hypergraph $H$ is obtained by replacing every vertex $v$ by $d$ vertices $v_1,\dots,v_d$. Its set of hyperedges is in bijection with that one of $H$, and this bijection is such that a hyperedge of $H$ contains $v$ if and only if its corresponding hyperedge in the inflation contains $v_1,\dots,v_d$; see Figure~\ref{fig: intro inflation}.

The structure of the paper is as follows. Section~2 presents notation, conventions, the fundamental notions concerning hypergraph associahedra, as well as a brief summary of Li's theory of wonderful compactifications. Section~3 contains the proof of Theorem~\ref{intro: main theorem mon}. In Section~4 we introduce the notion of inflation for hypergraphs, and briefly explain its combinatorial interpretation for associahedra and cyclohedra. We also present a streamlined account on compactifications $T_{d,n}^\cal{A}$, and conclude the section with a proof of Theorem~\ref{intro: main thm tdn}.

\subsection*{Acknowledgements.} The authors are very grateful to Patricio Gallardo and Jose Gonzalez for their constant support, and for feedback on the first draft of the paper. This paper builds on the first-named author's bachelor's thesis at Harvey Mudd College, supervised by the second-named author. We thank the College and its Mathematics Department for their supportive environment, where this work was conducted.

\section{Background}
\subsection{Notation and conventions}
Given $n\in\bb{Z}_{>0}$ fix the notation $[n]:=\{1,\dots,n\}$. Let $\{e_i\}_{i\in [n]}$ denote the standard basis of $\bb{R}^n$. Throughout we work with the realization of the $n$-dimensional simplex given by
\[
    \Delta_n := \operatorname{Conv}(\overline{e}_i\,\vert\, i\in [n+1])\subset\bb{R}^{n+1}/\bb{R}(1,\dots,1).
\]
Its inner normal fan, denoted $\Sigma_n$, is the collection of all cones of the form
\[
    \cone(\overline{e}_{i_1},\dots,\overline{e}_{i_k}\,\vert\, I=\{i_1,\dots,i_k\}\subsetneq [n]),
\]
for all $I\subsetneq[n]$; here $\cone(\emptyset)=\{0\}$. Finally, for $v=\sum_{i=1}^nv_ie_i\in\bb{R}^n$ and $I\subseteq [n]$, define
\[
    v_I := \sum_{i\in I} v_i\in\bb{R}.
\]

\subsection{Preliminaries on hypergraph associahedra}

Consider $n\in\bb{Z}_{>0}$. A \emph{hypergraph $H$ on $[n]$} is a subset $H\subseteq 2^{[n]}$ such that $\emptyset\not\in H$ and such that $\cup H = [n]$. The set $[n]$ is called the \emph{vertex set of $H$} and the elements of $H$ are called its \emph{hyperedges}. A hypergraph $H$ on $[n]$ is said to be:
\begin{itemize}
    \item \emph{Atomic}, if $\{i\}\in H$ for all $i\in [n]$.
    \item \emph{Saturated}, if given any $I,J\in H$ such that $I\cap J\neq\emptyset$, then $I\cup J\in H$.
\end{itemize}

A \emph{hypergraph partition} is a partition $\{H_1,\dots,H_n\}$ of $H$ such that $\{\cup H_1,\dots,\cup H_n\}$ is a partition of $\cup H$. Then, $H$ is said to be:

\begin{itemize}
    \item \emph{Connected}, if it has a unique hypergraph partition.
\end{itemize}

\begin{example}
    The hypergraph with vertex set $[4]$ and hyperedges $\{\{1,2\},\{3,4\}\}$ is saturated, but not atomic nor connected. The one with hyperedges $\{\{i\}\,\vert\,i\in[4]\}\cup\{\{1,2\},\{3,4\}\}$ is atomic and saturated, but not connected. Finally, the one with hyperedges $\{\{i\}\,\vert\,i\in[4]\}\cup\{\{1,2\},\{3,4\},\{1,2,3,4\}\}$ is atomic, saturated and connected.
\end{example}

\begin{definition}
    We use the capital letters A, S and C to describe a hypergraph $H$ that is atomic, saturated or connected. For example, an \emph{ASC-hypergraph} is one satisfying all these properties.
\end{definition}

With an AS-hypergraph as input, Do\v{s}en and Petri\'c define its correspond \emph{hypergraph associahedron} via their concepts of \emph{constructions} and \emph{constructs} \cite[Section~3]{dovsen2011hypergraph}. Here, however, we leverage the fact that hypergraph associahedra are known to be equivalent to another class of polytopes known as \emph{nestohedra}; see \cite[Sections~1~and~2]{petric2014stretching} for a dictionary between both theories. 

Let $H$ be an ASC-hypergraph with vertex set $[n]$. Then, it is possible to construct a complete fan $\Sigma(H)$ in $\bb{R}^{n+1}/\bb{R}(1,\dots,1)$ from it as follows.

\begin{construction}\label{bg: blow up construction}
    Start with the fan of the simplex $\Sigma_n\subset\bb{R}^{n+1}/\bb{R}(1,\dots,1)$, whose rays have been labeled by the singletons $i$ for $i\in [n]$. Given a nonsingleton hyperedge $I\in H$, define the cone $\sigma_I:=\cone(\overline{e}_i\,\vert\, i\in I)\in\Sigma_n$. Finally, define a total order of $H=\{I_1,\dots,I_k\}$ such that larger sets come before smaller sets. Then, $\Sigma(H)$ is the fan obtained from $\Sigma_n$ after performing an iterated star subdivision of the cones $\sigma_I$ as $I$ ranges through the hyperedges of $H$ in our choice of total order. 
\end{construction}

The fan $\Sigma(H)$ is independent of the choice of total order, by \cite[Theorem~4]{feichtner-muller}. Moreover, this fan is the inner normal fan of a polytope; see \cite[Remark~6.6]{postnikov-reiner-williams} and its preceding discussion. This justifies the following:

\begin{definition}
    Let $H$ be an ASC-hypergraph on $[n]$. Then, any polytope whose inner normal fan is $\Sigma(H)$ is called a \emph{hypergraph associahedron of $H$}. We denote any particular choice of it as $\mathscr{P}H$.
\end{definition}

It is worth noting that $\Sigma(H)$ is a complete unimodular fan in $\bb{R}^{n+1}/\bb{R}(1,\dots,1)$ with respect to the standard lattice $\bb{Z}^{n+1}/\bb{Z}(1,\dots,1)\subset\bb{R}^n/\bb{R}(1,\dots,1)$, by \cite[Theorem~5.1~and~Corollary~5.2]{zelevinsky-nested}.

\begin{figure}
\centering
    \includegraphics[width=0.7\textwidth]{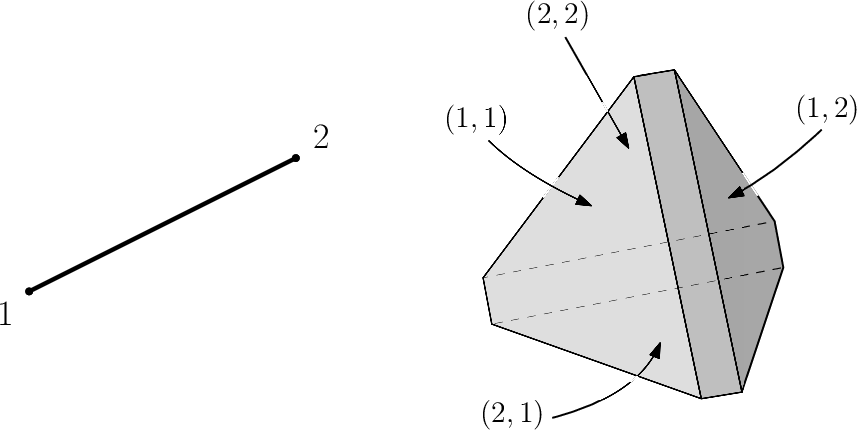}
    \caption{The $2$-inflation of the $1$-simplex. The vertex $i$ of the $1$-simplex gets ``inflated" to the facets $(i,1)$ and $(i,2)$ of its $2$-inflation, for $i=1,2$ (cf. Figure~\ref{fig: intro inflation}).}
    \label{fig:bg inflation}
\end{figure}

\begin{example}\label{bg: example P3}
    Consider the ASC-hypergraph with vertex set $[4]$ and having nonsingleton hyperedges $\{\{1,2\},\{3,4\},\{1,2,3,4\}\}$. Then, $\Sigma(H)$ is obtained by subdividing the cones $\cone(\overline{e}_1,\overline{e}_2)$ and $\cone(\overline{e}_3,\overline{e}_4)$ of the fan $\Sigma_3\subset\bb{R}^4/(1,1,1,1)$. This corresponds to truncating the $3$-simplex along two disjoint edges. See Figure~\ref{fig:bg inflation}.
\end{example}

\begin{remark}
    Hypergraph associahedra, also known as nestohedra, were originally introduced in \cite{dovsen2011hypergraph}, leveraging the theory of hypergraphs as an alternative to the use of combinatorial building sets. It follows directly from the definitions that an AS-hypergraph is the same as a combinatorial building set \cite[Definition~2.1]{zelevinsky-nested}, and an ASC-hypergraph is the same as a connected combinatorial building set. In this note we will adhere to the terminology of hypergraph associahedra, but use freely the equivalence of this concept and nestohedra to make use of some of the known results that appear in \cite{feichtner-strumfels,postnikov2009permutohedra, zelevinsky-nested}. We refer the reader to \cite{zelevinsky-nested,postnikov-reiner-williams} for more details on combinatorial building sets and nestohedra.
\end{remark}

\subsection{Graph associahedra} A prominent example of hypergraph associahedra are the so-called graph associahedra, originally defined by Carr and Devadoss \cite{carr-devadoss}. Important examples of graph associahedra include associahedra, cyclohedra and permutohedra. 

Given a connected graph $G$ with vertex set $[n]$, define the ASC-hypergraph $H(G)$ on $[n]$  as the hypergraph with hyperedges $\{I\subsetneq [n]\,\vert\, G|_I\text{ is connected}\}$. In the terminology of Carr-Devadoss, the elements of $H(G)$ are called \emph{tubes}.

\begin{definition}
    The \emph{graph associahedron of $G$} is the hypergraph associahedron of $H(G)$.
\end{definition}

Note that the process outlined in Construction~\ref{bg: blow up construction} is precisely that one used by Carr and Devadoss to define the graph associahedron of $G$.

\subsection{General geometric conventions} 

For the sake of brevity we will not delve into any details on toric geometry in this note, instead, we refer the reader to \cite{cox2011toric,fulton1993introduction}.

For simplicity we assume all our algebraic varieties are defined over $\bb{C}$. Given a positive integer $n$, we always consider the projective space $\bb{P}^n$ as the toric variety constructed from the fan $\Sigma_{n}\subset\bb{R}^{n+1}/\bb{R}(1,\dots,1)$. With these conventions, a nonzero vector $x=\sum_{i\in [n+1]} x_ie_i$ has homogeneous coordinates $[x_1:\dots:x_{n+1}]\in\bb{P}^n$. 

From our choice of fan, the torus-invariant divisors of $\bb{P}^n$ are precisely the hyperplanes $L_i:=V(x_i)$ for $i\in[n+1]$. The rest of the torus-invariant subvarieties are  
\[
    L_I := \bigcap_{i\in I} L_i.
\]
for every nonempty subset $I\subsetneq [n+1]$.

\begin{definition}\label{bg: def toric variety of hypergraph}
    Given a fan $\Sigma$ supported on $\bb{R}^{n+1}/\bb{R}(1,\dots,1)$, we denote its corresponding toric variety by $X(\Sigma)$. When $H$ is an SC-hypergraph (not necessarily atomic), we denote by $X(H)$ the toric variety corresponding to $\Sigma(H^{at})$, where $H^{at}$ is the ASC-hypergraph obtained from $H$ by adding all singletons as hyperedges.
\end{definition}

\begin{example}\label{bg: example P3 2}
    The toric variety $X(H)$ for $H$ as in Example~\ref{bg: example P3} is the blow-up of $\bb{P}^3$ along the two disjoint torus-invariant lines $L_{12}=L_1\cap L_2$ and $L_{34}=L_3\cap L_4\subset\bb{P}^3$.
\end{example}

\subsection{Wonderful compactifications}

To conclude the preliminaries we give a brief survey on Li's theory of wonderful compactifications \cite{li2009wonderful}. For us, Li's results play an essential role in the proof of Theorem~\ref{intro: main theorem mon}. Throughout this subsection we fix $Y$ to be a nonsingular variety over $\bb{C}$.

\begin{definition}{\cite[Section~5.1]{li2009wonderful}}
    Let $A,A_1,\dots,A_k,B$ be nonsingular subvarieties of $Y$. Then,
    \begin{itemize}
        \item The intersection of $A$ and $B$ is said to be \emph{clean} if the set-theoretic intersection $A\cap B$ is nonsingular, and we have the following relation between tangent spaces:
        \[
            T_{A\cap B,y}=T_{A,y}\cap T_{B,y},\text{ for all }y\in A\cap B.
        \]
        \item The intersection of $A_1,\dots,A_k$ is said to be \emph{transversal}, if either $k=1$, or for all $y\in Y$
        \[
            \operatorname{codim}\left( \bigcap_{i=1}^k T_{A_i,y},T_y\right)
            =
            \sum_{i=1}^k \operatorname{codim}(A_i,Y).
        \]   
    \end{itemize}
\end{definition}

\begin{definition}\cite[Definition~2.1]{li2009wonderful}
    An \emph{arrangement} of subvarieties of $Y$ is a finite collection $\cal{S}=\{S_i\}$ of properly contained nonsingular subvarieties $S_i\subsetneq Y$ such that:
    \begin{enumerate}
        \item[(i)] $S_i$ and $S_j$ intersect cleanly.
        \item[(ii)] $S_i\cap S_j$ is either equal to some $S_k$ or empty.
    \end{enumerate}
\end{definition}

In simple terms, the collection $\cal{S}$ in the previous definition is closed under intersections and these intersections are empty or clean.

\begin{definition}{\cite[Definition~2.2]{li2009wonderful}}
    Let $\cal{S}$ be an arrangement of subvarieties of $Y$. A subset $\cal{G}\subseteq\cal{S}$ is called a \emph{building set of $\cal{S}$} if, for all $S\in\cal{S}$, the minimal elements of $\{G\in\cal{G}\,\vert\, S\subseteq G\}$ intersect transversally, and their intersection is $S$.
\end{definition}

Note that in the previous definition the condition is automatically satisfied if $S\in\cal{G}$ by the definition of transversal intersection.

With these definitions, we can introduce the main objects that concern us.

\begin{definition}{\cite[Definition~2.2]{li2009wonderful}}
    A finite collection $\cal{G}$ of nonsingular subvarieties of $Y$ is called a \emph{building set} if the set $\cal{S}$ of all possible intersections of elements of $\cal{G}$ is an arrangement, and furthermore $\cal{G}$ is a building set of $\cal{S}$.
\end{definition}

The following is a direct consequence of the previous definitions. Indeed, the key is to note that the transversality assumption is satisfied when the intersection of two elements of a building set is also an element of the building set.

\begin{lemma}\label{bg: building set of itself}
    If $\cal{G}$ is an arrangement, then it is a building set of itself.
\end{lemma}

With these definitions at hand one can define the wonderful compactification of a building set.

\begin{definition}{\cite[Definition~1.1]{li2009wonderful}}\label{bg: def wonderful}
    Let $\cal{G}$ be a nonempty building set of subvarieties of $Y$, and $Y^\circ=Y\setminus \bigcup_{G\in\cal{G}}G$. Then, the closure of the image of the locally closed embedding
    \[
        Y^\circ\hookrightarrow\prod_{G\in\cal{G}}\operatorname{Bl}_G Y
    \]
    is called the \emph{wonderful compactification} of $\cal{G}$ on $Y$. It is denoted $Y_\cal{G}$.
\end{definition}

Prominent examples of building sets include Kapranov's construction of $\overline{M}_{0,n}$, discussed in Section~3, the Fulton-MacPherson configuration spaces \cite{fulton1994compactification}, and De Concini and Procesi's wonderful models of subspace arrangements \cite{de1995wonderful}.

For our purposes, we will use the fact that wonderful compactifications can be constructed as iterated blow-ups and, moreover, that the order in which these blow-ups are performed can often be changed without affecting the resulting variety.

\begin{proposition}{\cite[Theorem~1.3.(ii)]{li2009wonderful}}\label{bg: li thm 1.3}
    Let $\cal{G}=\{G_1,\dots,G_N\}$ be a nonempty building set of subvarieties of $Y$. Choose any total order on $\cal{G}$ such that $\cal{G}_r:=\{G_1,\dots,G_r\}$ is a building set for any $1\leq r\leq N$. Then, 
    \[
        Y_\cal{G} = \operatorname{Bl}_{\widetilde{G}_N}\cdots \operatorname{Bl}_{\widetilde{G}_2}\operatorname{Bl}_{G_1}Y,
    \]
    where $\widetilde{G}$ is the strict transform of $G$, and each blow-up is performed along a nonsingular subvariety.
\end{proposition}

To conclude these remarks let us further explore some consequences of the previous result in the context of linear subspaces of projective space, which is the case that we will be using.

\begin{lemma}\label{bg: lemma linear subspaces closed under intersection implies building set}
    Let $Y=\bb{P}^n$ and $\cal{G}=\{L_1,\dots,L_N\}$ be a collection of linear subspaces that is closed under intersection. Then, $\cal{G}$ is an arrangement. In particular, by Lemma~\ref{bg: building set of itself}, it is a building set of itself.
\end{lemma}

\begin{proof}
    Since $\cal{G}$ is closed under intersection, all that remains to show is that its elements intersect cleanly. However, this is a direct consequence of the fact that the $L_i$ are linear subspaces. Indeed, by identifying $L_i$ with its embedded tangent space one obtains $T_{L_i,y}\cap T_{L_j,y} = L_i\cap L_j = T_{L_i\cap L_j,y}$.
\end{proof} 

In light of Proposition~\ref{bg: li thm 1.3} it is important to determine some blow-up orders that give rise to the wonderful compactification.

\begin{lemma}\label{bg: total orders giving wonderful compactification}
    Let $Y=\bb{P}^n$ and $\cal{G}=\{L_1,\dots,L_N\}$ be a collection of linear subspaces that is closed under intersection. Then, 
    \[
        Y_\cal{G} = \operatorname{Bl}_{\widetilde{G}_N}\cdots \operatorname{Bl}_{\widetilde{G}_2}\operatorname{Bl}_{G_1}Y
    \]
    for any total order of $\cal{G}$ such that $\cal{G}_r=\{L_1,\dots,L_r\}$ is closed under intersection for all $1\leq r\leq N$.
\end{lemma}

\begin{proof}
    If $\cal{G}_r$ is closed under intersection, then, by Lemma~\ref{bg: lemma linear subspaces closed under intersection implies building set}, $\cal{G}_r$ is a building set. The result then follows from Proposition~\ref{bg: li thm 1.3}.
\end{proof}

\section{Hypergraph associahedra and Hassett spaces}
\subsection{Introduction} 

Consider the moduli problem of parametrizing $n$ distinct labeled points in $\bb{P}^1$ up to change of coordinates, this is, up to the action of $\operatorname{PGL}_2$ (equivalently, up to M\"obius transformations). Any configuration of distinct points $p_1=[u_1:v_1],\dots,p_n=[u_n,v_n]$ is equivalent to one of the form $[1:y_1],\dots,[1:y_{n-2}],[1:0],[0:1]$, where the $y_i$ are still only defined up to scaling. In other words, the space parametrizing these equivalence classes is
\begin{align*}
    M_{0,n} &= (\bb{P}^1)^n-\bigcup \{p_i=p_j\}/\operatorname{PGL_2}\\
    &= \{(y_1,\dots,y_{n-2})\in\bb{C}^{n-2}\,\vert\,y_i\neq 0,\text{ and }y_i\neq y_j\}/\bb{C}^\times.
\end{align*}

This space is manifestly not compact. However, Kapranov \cite[Section~4.2]{Kap93Chow} uses it to construct the Grothendieck-Knudsen compactification $\overline{M}_{0,n}$ as follows. First, embed the previous realization $M_{0,n}\subset\bb{P}^{n-3}=\bb{C}^{n-3}-\{0\}/\bb{C}^\times$ using the coordinates above. Then, Kapranov proves that $\overline{M}_{0,n}$ is the iterated blow-up of this $\bb{P}^{n-3}$ along all linear subspaces generated by the points $q_0=[1:\dots:1]$, and $q_1=[1:0:\cdots:0],\cdots,q_{n-2}=[0:\cdots:0:1]$ in non-decreasing order of dimension (blow-up the points first, then the strict transform of the lines joining any two of them, then strict transform of the planes through any three of them, etc.). 

The open subset $M_{0,n}\subset\overline{M}_{0,n}$ is often referred to as the \emph{interior} of $\overline{M}_{0,n}$, while its complement is often called the \emph{boundary} of $\overline{M}_{0,n}$. The key idea in the construction of $\overline{M}_{0,n}$ is to enlarge the moduli problem to that one of parametrizing $n$ points in \emph{stable} curves of genus $0$, up to componentwise change of coordinates.

\begin{definition}\label{mon: stable curve}
    An \emph{$n$-pointed stable rational curve} is a tuple $(C,p_1,\dots,p_n)$ such that
    \begin{enumerate}
        \item $C$ is a connected curve of arithmetic genus zero. The singularities of $C$ can only be nodes.
        \item The points $p_1,\dots,p_n\in C$ and the nodes of $C$ are all distinct.
        \item Every irreducible component of $C$ contains at least $3$ special points, where a special point is either a marked point or a node.
    \end{enumerate}
\end{definition}

The boundary of $\overline{M}_{0,n}$ parametrizes $n$-pointed stable rational curves with at least two irreducible components, up to componentwise automorphisms. More concretely, given any set $I\subseteq [n-1]$ such that $|I|\geq 2$, there is a divisor $D(I\,\vert\, I^C)\subset\overline{M}_{0,n}$ whose general point parametrizes an $n$-pointed stable rational curve with two components, one containing the points labeled by $I$, and the other one those labeled by $I^C$, in a such a way that the $(n-1)$th and $n$th points are in different components.

\subsection{Hassett compactifications} 

The compactification $\overline{M}_{0,n}$ is a particular example of a general class of compactifications of $M_{0,n}$ defined by Hassett \cite{hassett2003moduli}. Hassett's approach is to endow each marked point with a weight and define the notion of stability in this setting. 

\begin{definition}
    The \emph{domain of admissible weights} is the set
    \[
        \mathcal{D}_{0,n} = \{(a_1,\dots,a_n)\in\bb{Q}^n\,\vert\, 0<a_i\leq 1,\, a_{[n]}>2\}.
    \]
    Its elements are called \emph{collections of weight data}.
\end{definition}

\begin{definition}
    Given $\cal{A}\in\mathcal{D}_{0,n}$, we say an $n$-pointed rational curve is $\cal{A}$-stable if it satisfies properties (1) and (2) from Definition~\ref{mon: stable curve}, and also satisfies the property
    \begin{enumerate}
        \item[(3')] In every irreducible component of $C$, the total weight of its marked points together with its number of nodes must be strictly larger than $2$.
    \end{enumerate}
\end{definition}

For every $\cal{A}\in\cal{D}_{0,n}$, Hassett defines a compactification $\overline{M}_{0,\cal{A}}$ of $M_{0,n}$ parametrizing the $n$ distinct pointed rational $\cal{A}$-stable curves up to automorphism. 

\subsection{Reduction maps and chamber decompositions} 

Given two collections of weight data $\cal{A}=(a_1,\dots,a_n)$ and $\cal{B}=(b_1,\dots,b_n)\in\mathcal{D}_{0,n}$, write $\cal{A}\geq\cal{B}$ if $a_i\geq b_i$ for all $i\in [n]$.

\begin{proposition}{\cite[Theorem~4.1]{hassett2003moduli}}\label{mon: hassett reduction maps}
    Consider two collections of weight data $\cal{A},\cal{B}\in\mathcal{D}_{0,n}$ such that $\cal{A}\geq\cal{B}$. Then, there exists a natural reduction morphism
    \[
        \rho_{\cal{B},\cal{A}}:\overline{M}_{0,\cal{A}}\to\overline{M}_{0,\cal{B}}.
    \] 
\end{proposition}

\begin{proposition}{\cite[Proposition~4.4]{hassett2003moduli}} \label{mon: reduction composition}
If $\cal{A},\cal{B},\cal{C}\in\mathcal{D}_{0,n}$ are such that the reduction maps $\rho_{\cal{B},\cal{A}},\rho_{\cal{C},\cal{B}}$ and $\rho_{\cal{C},\cal{A}}$ are defined, then 
\[
    \rho_{\cal{C},\cal{A}} = \rho_{\cal{C},\cal{B}}\circ\rho_{\cal{B},\cal{A}}.
\]
\end{proposition}

A \emph{chamber decomposition} of the domain of admissible weights $\cal{D}_{0,n}$ is a finite collection $\cal{W}$ of hyperplanes or \emph{walls} $w_I\subset\cal{D}_{0,n}$. For a set $I\subsetneq [n]$, the wall $w_I$ is defined by the equation $a_I=1$, and the \emph{open chambers} are the connected components of $\cal{D}_{0,n}\setminus \cup w_I$. Hassett considers two chamber decompositions: the \emph{fine chamber decomposition}, consisting of all walls $w_I$ where $2\leq |I|\leq n-2$, and the \emph{coarse chamber decomposition}, defined by all walls $w_I$ where $2<|I|<n-2$. The key results are that $\overline{M}_{0,\cal{A}}\cong\overline{M}_{0,\cal{B}}$ if $\cal{A}$ and $\cal{B}$ are in the same coarse chamber, while these two spaces have the same universal family if $\cal{A}$ and $\cal{B}$ are in the same fine chamber \cite[Proposition~5.1]{hassett2003moduli}. 

\begin{definition}\label{mon: def geq_f}
    Given a collection of admissible weights $\cal{A}\in\cal{D}_{0,n}$ we denote its coarse (resp. fine) chamber by $Ch_c(\cal{A})$ (resp. $Ch_f(\cal{A})$). Given any other $\cal{B}\in\cal{D}_{0,n}$, we write $\cal{A}\geq_c \cal{B}$ if there are collections of weight data $\cal{A}'$ and $\cal{B}'$ in the same coarse chambers as $\cal{A}$ and $\cal{B}$, respectively, such that $\cal{A}'\geq\cal{B}'$.
\end{definition}

\begin{definition}
    Let $\cal{A},\cal{B}\in\cal{D}_{0,n}$. We say that there is a \emph{simple wall crossing from $\cal{A}$ to $\cal{B}$ along the wall $w_I$} if $\cal{A}\geq\cal{B}$, and there is a \emph{unique} $I\subset [n]$, $2\leq|I|\leq n-1$ such that $a_I>1$, but $b_I\leq 1$.
\end{definition}

The following two results will be essential for the proof of Theorem~\ref{intro: main theorem mon}.

\begin{lemma}\cite[Lemma~4.2]{alexeev2008moduli}\label{mon: lemma alexeev}
    Given any two collections of weight data $\cal{A}\geq\cal{B}$ there exist collections of weight data $\cal{A}'\geq\cal{B}'$ in the same coarse (fine) chambers as $\cal{A}$ and $\cal{B}$, respectively, such that and the straight line from $\cal{A}'$ to $\cal{B}'$ goes through a sequence of simple wall crossings.
\end{lemma}

\begin{proposition}{\cite[Remark~4.6 and Corollary~4.7]{hassett2003moduli}}\label{mon: hassett blowup}
        Suppose that $\cal{A}\geq\cal{B}$ and there is a simple wall crossing from $\cal{B}$ to $\cal{A}$ along the wall $w_I$. If $b_J\leq 1$ for every proper subset $J\subsetneq I$, then $\rho_{\cal{B},\cal{A}}$ is the blow-up of $\overline{M}_{0,\cal{B}}$ along the image of the boundary divisor $D(I\,\vert\, I^C)\subset\overline{M}_{0,\cal{A}}$ in $\overline{M}_{0,\cal{B}}$. In particular, if $|I|=2$, then $\rho_{\cal{B},\cal{A}}$ is a natural isomorphism $\overline{M}_{0,\cal{A}}\cong\overline{M}_{0,\cal{B}}$.
\end{proposition}

\subsection{The Losev-Manin and projective space as Hassett spaces}

The wall and chamber decompositions of $\cal{D}_{0,n}$ is incredibly complex. For example, the number of fine chambers grows proportionally with $2^{n^2-n\log_2 n+O(n)}$, by \cite[Theorem~1.3]{ascher2020enumerating}. Moreover, there is no guarantee that two different coarse chambers will have non-isomorphic coarse moduli spaces. Therefore, it is important to be precise on which realizations of the moduli spaces involved we are considering.

Let us fix $\overline{M}_{0,n}$ to be the Hassett space corresponding to $(1,\dots,1)\in\cal{D}_{0,n}$. Then, by \cite[Section~6.1]{hassett2003moduli}, Kapranov's blow-up construction $\overline{M}_{0,n}\to\bb{P}^{n-3}$, as outlined at the beginning of this section, is realized by the reduction map $\overline{M}_{0,n}\to\bb{P}^{n-3}=\overline{M}_{0,\cal{A}_\bb{P}}$, where $\cal{A}_\bb{P}=(1/(n-2),\dots,1/(n-2),1)$. However, considering only weights $\cal{A}$ such that $\cal{A}\geq\cal{A}_\bb{P}$ in Theorem~\ref{mon: main theorem} is an unnecessary restriction.

Note that when $n=4$ the reduction map $\overline{M}_{0,n}\to\bb{P}^1$ is an isomorphism, as every blow-up of $\bb{P}^1$ along a finite set of points leaves the space unchanged. Therefore, from hereafter many of our results will explicitly consider $n\geq 5$.

\begin{definition}\label{mon: def pn}
    Define $\cal{A}_{\bb{P}}:=(1/(n-2),\dots,1/(n-2),1)\in\cal{D}_{0,n}$. Following \cite[Section~6.1]{hassett2003moduli}, from hereafter we consider $\bb{P}^{n-3}$ as a Hassett space corresponding to any collection of weight data defined by the inequalities
    \begin{align}\label{mon: projective space ineqs}
    \begin{split}
        &a_{[n-1]\setminus\{i\}}\leq 1,\text{ for all }i\in[n-1],\\
        &a_{[n-1]}>1,
    \end{split}
    \end{align}
    If $n\geq 5$, the interior of this set is a coarse chamber, which we denote by $Ch_c(\cal{A}_{\bb{P}})$. 
\end{definition}

\begin{remark}
    Strictly speaking $\cal{A}_{\bb{P}}$ lies on the boundary of $Ch_c(\cal{A}_{\bb{P}})$. However, by \cite[Section~6.1]{hassett2003moduli}, every collection of weight data satisfying the inequalities of Definition~\ref{mon: def pn} defines isomorphic coarse moduli spaces.
\end{remark}

In the same spirit, let us define the realizations of the Losev-Manin space we consider here.

\begin{definition}\label{mon: def LM}
    Define $\cal{A}_{LM}:=(1/(n-2),\dots,1/(n-2),1,1)\in\cal{D}_{0,n}$. Following \cite[Section~6.4]{hassett2003moduli}, from hereafter we consider $\overline{M}_{0,n}^{LM}$ as a Hassett space corresponding to any collection of weight data defined by the inequalities
    \begin{align}\label{mon: LM ineqs}
    \begin{split}
        &a_n+a_i>1,\, a_{n-1}+a_i>1,\text{ for all }i\in[n],\\
        &a_{[n-2]}\leq 1.
    \end{split}
    \end{align}
    If $n\geq 5$, the set defined by the last inequality is a coarse chamber, which by abuse of notation we denote by $Ch_c(\cal{A}_{LM})$. 
\end{definition}

\begin{remark}
    Strictly speaking $\cal{A}_{LM}$ lies on the boundary of $Ch_c(\cal{A}_{LM})$. However, by \cite[Section~6.4]{hassett2003moduli}, every collection of weight data satisfying the inequalities of Definition~\ref{mon: def LM} defines isomorphic coarse moduli spaces.
\end{remark}

As observed in \cite[Section~6.4]{hassett2003moduli}, the composition of reduction maps $\overline{M}_{0,n}\to\overline{M}_{0,n}^{LM}\to\bb{P}^{n-3}$ induced by $(1,\dots,1)\geq\cal{A}_{LM}\geq\cal{A}_{\bb{P}^{n-3}}$ gives a construction of $\overline{M}_{0,n}^{LM}$ as an iterated blow-up along all linear subspaces generated by the points $q_1=[1:0:\cdots:0],\cdots,q_{n-2}=[0:\cdots:0:1]$ in non-decreasing order of dimension (blow-up the points first, then the strict transform of the lines joining any two of them, then strict transform of the planes through any three of them, etc.). In particular, $\overline{M}_{0,n}^{LM}$ is a toric variety with the standard toric structure introduced in the Background section. From Construction~\ref{bg: blow up construction} one readily obtains that its corresponding fan is $\Sigma(H)$ for the ASC-hypergraph $H$ having vertex set $[n-2]$ and hyperedge set $2^{[n-2]}$.

\subsection{Proof of Theorem~\ref{intro: main theorem mon}}

Ferreira da Rosa, Jensen and Ranganathan \cite{da2016toric} gave a full classification of which graph associahedra give rise to a toric variety that is isomorphic to a Hassett space and, conversely, of which Hassett spaces are isomorphic to the toric variety of a graph associahedron. In this section we generalize their latter result by showing that every Hassett space ``between" projective space and the Losev-Manin compactification is the toric variety of a hypergraph associahedron.

Every $d$-dimensional hypergraph associahedron is constructed as an iterated truncation of the $d$-simplex, so that the corresponding toric variety is the iterated blow-up of $\bb{P}^{d}$. Consequently, we only consider compactifications that admit a reduction map to projective space. In particular, this means that our results fail to classify all possible toric compactifications, as evidenced by Keel's construction of $\overline{M}_{0,n}$ as a an iterated blow-up of the Hassett space $(\bb{P}^1)^{n-3}$; see \cite[Section~6.3]{hassett2003moduli}.

\begin{theorem}\label{mon: main theorem}[Theorem~\ref{intro: main theorem mon}]
    Let $n\geq 4$ and $\cal{A}\in\cal{D}_{0,n}$ be such that $\cal{A}_{LM}\geq_c \cal{A}$ and $\cal{A}\geq_c \cal{A}_\bb{P}$. Then, $\overline{M}_{0,\cal{A}}\cong X(H_\cal{A})$ for the SC-hypergraph $H_\cal{A}$ in the following definition.
\end{theorem}

\begin{definition}\label{mon: def HA}
    Let $n\geq 4$ and $\cal{A}\in\cal{D}_{0,n}$ be such that $\cal{A}_{LM}\geq_c \cal{A}$ and $\cal{A}\geq_c \cal{A}_\bb{P}$. Define the $H_\cal{A}$ to be the SC-hypergraph with vertex set $[n-2]$ and hyperedges
    \[
        H_\cal{A}:=\left\{I\,\vert\,I\subseteq[n-2],\,|I|\geq 1\text{ and }a_I+a_{n-1}>1\right\}.
    \]  
    It is saturated because if $a_I+a_{n-1}>1$ and $a_J+a_{n-1}>1$, then $a_{I\cup J}+a_{n-1}>1$. Connectedness follows because $a_{[n-1]}>1$ for every element in $\cal{D}_{0,n}$, since $a_n\leq 1$ and $a_{[n]}>2$.
\end{definition}

\begin{remark}
    In light of the theorem, it is natural to ask whether every SC-hypergraph gives rise to a toric compactification of $M_{0,n}$ as in the statement. The answer to this question is negative, as for a general SC-hypergraph on $[n]$ it is possible that a hyperedge of $H$ is contained in a proper subset $S\subseteq [n]$ that is not a hyperedge, which is impossible in our context. 
\end{remark}

Let us present some essential results for the proof of the main theorem. 

\begin{lemma}\label{mon: lemma chambers pn}
    Suppose $n\geq 5$, and let $\cal{A}=(a_1,\dots,a_n)\in\cal{D}_{0,n}$. Then,
    \begin{enumerate}
        \item[1)] we always have that $a_{[n-1]}>1$.
        \item[2)] If $\cal{A}\geq_c \cal{A}_\bb{P}$, then $a_I+a_n>1$ for all $I\subseteq[n-1]$, $|I|\geq 2$.
        \item[3)] If $\cal{A}_{LM}\geq_c \cal{A}$, then $a_I\leq 1$ for all $I\subseteq [n-2]$.
    \end{enumerate}
    In particular, all these inequalities hold true in the whole coarse chamber $Ch_c(\cal{A})$ in their respective cases.
\end{lemma}
\begin{proof}
    The proof of (1) is an immediate consequence of the fact that $a_n\leq 1$ and $a_{[n]}>2$. We only present the proof of (2), the proof of (3) being similar.
    
    Let $\cal{A}'\in Ch_c(\cal{A})$ and $\cal{B}'=(b_1,\dots,b_n)\in Ch_c(\cal{A}_{\bb{P}})$ be such that $\cal{A}'\geq\cal{B}'$. Then, from Definition~\ref{mon: def pn}, we see that $b_{[n-1]\setminus\{i\}}\leq 1$ for all $i\in[n-1]$, as $|[n-1]\setminus\{i\}|\geq 3$. It follows that $b_i+b_n>1$, since $ b_{[n]}>2$. In particular, this implies that $b_I+b_n>1$ for all $I\subseteq [n-1]$, $|I|\geq 1$, and the same inequalities hold for $\cal{A}'$. Since $\cal{A}'$ and $\cal{A}$ are in the same coarse chamber, their components satisfy the same inequalities provided that these inequalities have three or more terms. Therefore, we conclude that $a_I+a_n>1$ for all $I\subseteq [n-1]$, $|I|\geq 2$.
\end{proof}

The following lemma is a key component in the proof of Theorem~\ref{mon: main theorem}.

\begin{lemma}\label{mon: lemma n-1 walls crossing}
    Let $n\geq 5$, and $\cal{A}\in\cal{D}_{0,n}$ be such that $\cal{A}_{LM}\geq_c\cal{A}$ and $\cal{A}\gneq_c \cal{A}_\bb{P}$. Then, the straight line from $Ch_c(\cal{A}_{\bb{P}})$ to $Ch_c(\cal{A})$ given by Lemma~\ref{mon: lemma alexeev} may only cross simple walls $w_J$ such that
    \begin{align*}
    J = \left\{
        \begin{aligned}
            &I\cup\{n-1\},\text{ for }I\subsetneq[n-2]\text{ with }|I|\geq 1;\text{ or}\\
            &\{i,n\}\text{ for }i\in[n-2]. 
        \end{aligned}
    \right.
    \end{align*}
\end{lemma}
\begin{proof}
    In order, the clauses of Lemma~\ref{mon: lemma chambers pn} imply that both $\cal{A}_{\bb{P}}$ and $\cal{A}$ lie on:
    \begin{enumerate}
        \item[1)] the ``$>1$" side of the walls $w_{I\cup\{n\}}$, for all $I\subseteq[n-1]$, $|I|\geq 2$;
        \item[2)] the ``$>1$" side of the wall $w_{[n-1]}$;
        \item[3)] the ``$\leq 1$" side of the walls $w_I$, for all $I\subseteq [n-2]$.
    \end{enumerate}
    Item (3) implies that if the wall $w_J$ is crossed from $\cal{A}_\bb{P}$ to $\cal{A}$, then $J$ must contain either $n-1$ or $n$. If $n\in J$, then Item (1) forbids $|J|\geq 3$, so $J=\{i,n\}$ for some $i\in[n-1]$. First note that the wall $w_J$ for $J=\{n-1,n\}$ is never crossed, as $Ch_c(\cal{A}_\bb{P})$ lies on its ``$>1$" side of it, as every element in this chamber satisfies $a_{[n-2]}\leq 1$ by Item (3). Finally, if $n-1\in J$ and $J\neq \{n-1,n\}$, then it follows from (2) that $J\subsetneq [n-1]$, so that $J=I\cup\{n-1\}$ for some $I\subsetneq[n-2]$, $|I|\geq 1$.
\end{proof}

\begin{proof}[Proof of Theorem~\ref{mon: main theorem}]
    If $n=4$, then $\overline{M}_{0,n}^{LM}=\overline{M}_{0,\cal{A}_{\bb{P}}}=\bb{P}^1$ and the result is immediate. Let us then suppose that $n\geq 5$. First consider the case in which $\cal{A}\in Ch_c(\cal{A}_{\bb{P}})$, so that $\overline{M}_{0,\cal{A}}\cong\bb{P}^{n-3}$. Then, $H_\cal{A}=\{\{i\}\,\vert\,i\in[n-2]\}\cup\{[n-2]\}$. Indeed, $I=[n-2]$ is in $H_\cal{A}$ by connectedness. Similarly, every singleton is in $H_\cal{A}$ by the argument given in the proof of Lemma~\ref{mon: lemma chambers pn}. Lastly, note that if $I\subsetneq [n-2]$ with $|I|\geq 2$, then $I\cup\{n-1\}\subsetneq [n-1]\setminus\{i\}$ for some $i\in [n-2]\setminus I$, so that $a_I+a_{n-1}\leq a_{[n-1]\setminus\{i\}} \leq 1$. It follows from Construction~\ref{bg: blow up construction} that $X(H_\cal{A})=\bb{P}^{n-3}$. Let us now prove the remaining cases.

    By Proposition~\ref{mon: hassett blowup} and Lemma~\ref{mon: lemma alexeev}, one can conclude that the reduction map $\rho:\overline{M}_{0,\cal{A}}\to\bb{P}^{n-3}$ is a composition of blow-up maps (when crossing walls in the coarse chamber decomposition) and natural isomorphisms (when crossing walls in the fine chamber decomposition that are not in the coarse one), one per simple wall crossing in the straight line connecting $Ch_c(\cal{A})$ and $Ch_c(\cal{A}_{\bb{P}})$. We show that these iterated blow-ups yield the same variety as the blow-up construction of $X(H_\cal{A})$ given in Construction~\ref{bg: blow up construction}.
    
    By Lemma~\ref{mon: lemma n-1 walls crossing}, every wall crossed along this path is of the form $w_{I\cup\{n-1\}}$, where $I\subsetneq [n-2]$, $|I|\geq 1$, or $w_{\{i,n\}}$ for some $i\in[n-1]$. Crossing the walls of the first form with $|I|=1$, or those of second form, defines a natural isomorphism between the coarse moduli spaces on both sides of the wall, so they are inconsequential to our analysis. Define the SC-hypergraph
    \[
        H_{\cal{I}}=\{ I\,\vert\, I\subsetneq [n-2], |I|\geq 1,\text{ and } w_{I\cup\{n-1\}}\text{ was crossed}\}\cup\{[n-2]\}.
    \]
    We claim that $H_\cal{I}=H_\cal{A}$. First note that $[n-2]$ is a hyperedge for both of them. On the other hand, note that $\cal{A}_\bb{P}$ lies on the ``$\leq 1$" side of the wall $w_{I\cup\{n-1\}}$ for all $I\subsetneq [n-2], |I|\geq 1$, because the same is true for the wall $w_{[n-1]\setminus\{i\}}$ for all $i\in[n-1]$, and $I\cup\{n-1\}\subseteq [n-1]\setminus\{i\}$ for some $i$. Then, for all $I\subsetneq [n-2]$ with $|I|\geq 1$ the wall $w_{I\cup\{n-1\}}$ is crossed if and only if $a_{I\cup\{n-1\}}>1$. Therefore, the defining conditions for both hypergraphs are equivalent.
    
    It only remains to show that the iterated blow-up map $\rho:\overline{M}_{0,\cal{A}}\to\bb{P}^{n-3}$ yields the same variety as the SC-hypergraph $H_\cal{A}$ via the blow-up construction $\rho':X(H_\cal{A})\to\bb{P}^{n-3}$. 

    Cnsider $J=I\cup\{n-1\}$ with $I\subsetneq[n-2]$, $|I|\geq 2$, and such that the wall $w_J$ was crossed. Then, by \cite[Section~6.2]{hassett2003moduli}, under the reduction map $\overline{M}_{0,\cal{A}}\to\bb{P}^{n-3}$, the divisor $D(I\cup\{n-1\}\,\vert\,I^C\cup\{n\})\subset\overline{M}_{0,\cal{A}}$ is the (iterated) proper transform of the linear subspace generated by the coordinate points indexed by $I^C$. In turn, this linear subspace corresponds to the cone $\operatorname{Cone}(\overline{e}_i\,\vert\, i\in I)$ in the fan $\Sigma_{n-3}$ of $\bb{P}^{n-3}$ under the toric orbit-cone correspondence.
    
    By the previous discussion both maps are iterated blow-ups of the same torus-invariant subvarieties of $\bb{P}^{n-3}$. Let us proceed to show that the order in which these blow-ups are performed is inconsequential. 
    
    Since blowing-up codimension-$1$ subvarieties is an isomorphism, in both cases the blow-ups are performed along the collection of (strict transforms of the) linear subspaces in the collection
    \[
        \cal{G}=\left\{L_I:=\bigcap_{i\in I}L_i\text{ for all }I\in H_\cal{I},\, 1<|I|<n-2\right\}.
    \]
    On the one hand, any total order of non-decreasing dimension on $\cal{G}$ can be used to construct $X(H_\cal{A})$, by \cite[Theorem~4]{feichtner-yuzvinsky}. On the other hand, as discussed at the beginning of the proof, $\overline{M}_{0,\cal{A}}$ is the iterated blow-up of $\bb{P}^{n-3}$ along this collection of subvarieties in the total order prescribed by the straight line from Lemma~\ref{mon: lemma alexeev}, starting from a wall adjacent to $Ch_c(\cal{A}_{\bb{P}})$ and finishing with the wall adjacent to $Ch_c(\cal{A})$. We will use Li's results on wonderful compactifications discussed in the background section to prove that the resulting varieties are isomorphic. 

    The first step is to see that, by Lemma~\ref{bg: lemma linear subspaces closed under intersection implies building set}, $\cal{G}$ is a (geometric) building set whose induced arrangement is $\cal{G}$ itself. Indeed, $\cal{G}$ is closed under intersection because if $I,J\subsetneq[n-2]$ and $1<|I|,|J|<n-2$, then $L_I\cap L_J=L_{I\cup J}$, which is empty if $I\cup J=[n-2]$. 

    Since $\cal{G}$ is a building set consisting of linear subspaces, recall that Lemma~\ref{bg: total orders giving wonderful compactification} says that if one reorders $\cal{G}$ in such a way that $L_{I_1},\dots,L_{I_r}$ is closed under intersection for all $1\leq r\leq N$, then the variety obtained by blowing up $\bb{P}^{n-3}$ along the elements of $\cal{G}$ in this order is isomorphic to the wonderful compactification $(\bb{P}^{n-3})_\cal{G}$. We claim that both blow-up orders in the constructions of $\overline{M}_{0,\cal{A}}$ and $X(H_\cal{A})$ satisfy this property, so both are isomorphic to $(\bb{P}^{n-3})_\cal{G}$.

    First note that any total order of non-decreasing dimension used to construct $X(H_\cal{A})$ satisfies this property. For the case of $\overline{M}_{0,\cal{A}}$, observe that if $L_I,L_J\in\cal{G}_r$, then $L_{I\cup J}=L_I\cap L_J$ appears before them in $\cal{G}_r$, since it is impossible to cross either of the walls $w_I$ and $w_J$ before crossing the wall $w_{I\cup J}$. Therefore, both constructions yield isomorphic varieties $(\bb{P}^{n-3})_\cal{G}\cong\overline{M}_{0,\cal{A}}\cong X(H_\cal{A})$.
\end{proof}

\section{Inflation and the moduli space of labeled points in affine space up to translation and scaling}\label{sec: tdn}

In this section we prove Theorem~\ref{intro: main thm tdn}. In order to do this, we introduce the notion of inflation of a hypergraph associahedron. Afterwards we present an overview of the moduli spaces under consideration.

\subsection{Inflating hypergraph associahedra}

The definition of inflation arose in our studies of the moduli space $T_{d,n}^{LM}$ that is introduced later in this section. However, the construction is general enough to be of independent interest.

\begin{definition}\label{tdn: def inflation}
    Consider a tuple of positive integers $\mbf{d}=(d_1,\dots,d_n)$, and let $H$ be an SC-hypergraph on $[n]$ (not necessarily atomic). Then, the \emph{$\mbf{d}$-inflation of $H$} is the ASC-hypergraph $\operatorname{Inf}_\mbf{d}(H)$ with vertex set $\coprod [d_i]$ and hyperedge set
    \[
        \operatorname{Inf}_\mathbf{d}(H) = \left\{\{j\}\,\vert\,j\in\coprod[d_i]\right\}\cup\left\{ \coprod_{i\in I}[d_i]\,\vert\, I\in H\right\}.
    \]
    If $\mbf{d}=(d,\dots,d)$ is constant, then we write $\operatorname{Inf}_d(H)$ and refer to it as the $d$-inflation of $H$.
\end{definition}   

\begin{example}
    Consider the hypergraph $H(K_n)$ corresponding to the complete graph $G=K_n$ on $n$ vertices. Then, its corresponding hypergraph associahedron is the $(n-1)$-dimensional permutohedron \cite{carr-devadoss}. Let us compute its $d$-inflation, which by abuse of notation we denote as $\operatorname{Inf}_d(K_n)$. By definition, its vertex set is the disjoint union of $n$ copies of $[d]$, this is, $\coprod [d]$. Its hyperedge set is
    \[
        \operatorname{Inf}_d(K_n) = \left\{\{i\}\,\vert\,i\in \coprod [d]\right\}\cup\left\{ \coprod_{i\in I}[d]\,\vert\, I\subsetneq [n],\, I\neq\emptyset\right\}.
    \]
    Figure~\ref{fig: intro inflation} in the Introduction shows the $2$-inflation of $K_2$.
\end{example}

\begin{example}\label{inf: example inflation}
    If $\cal{P}G$ is an associahedron or cyclohedron with vertex set $[n]$, then the $d$-inflation has a straightforward combinatorial interpretation. Recall that the associahedron is the graph associahedron of the path graph, while the cyclohedron is that one of the cycle graph. In both cases a tube of $G$ corresponds to a parenthesization of the vertices of their respective graphs. The only difference between both cases is that the parenthesization in the path graph is restricted by its endpoints, while for the cycle graph it is not. Then, the faces of the $d$-inflation of the graph associahedron $\cal{P}(\operatorname{Inf}_d(H(G))$ correspond to parenthesizations of a path, respectively cycle, graph with $dn$ vertices, in such a way that the parenthesizations occur in fixed blocks of size $d$. These blocks cannot be split by any set of parenthesis, however. 
\end{example}

\subsection{Toric compactifications of the moduli space of points in affine space arising through the operation of inflation} 

Chen, Gibney and Krashen \cite{Chen-Gibney-Krashen} compactified the moduli problem of parametrizing $n$ distinct labeled points in $d$-dimensional affine space up to translation and scaling. They denoted their moduli space $T_{d,n}$. Following up on their work, Gallardo and Routis \cite{gallardo2017wonderful} constructed many more compactifications by attaching weights to the points and allowing collisions between them depending on their weights. Let us briefly recall some of the main results concerning Gallardo and Routis' weighted compactifications. 

\subsection{Gallardo and Routis' compactifications} Consider the problem of parametrizing $n$ distinct labeled points in the affine space $\bb{A}_\bb{C}^d$ up to translation and scaling. More precisely, we are interested in studying the geometry of the orbit space of the diagonal action of $\CC^\times\ltimes\CC^d$ on $(\bb{A}_\bb{C}^d)^n$ defined in each component as $(t,\vec{s})\cdot p = t\cdot p + \vec{s}$. Configurations of distinct labeled points up to this action are parametrized by an open subset $T_{d,n}^{\circ}\subseteq\bb{P}^{d(n-1)-1}$ consisting of the points
\[
    [x_{11}: \cdots : x_{1d} : x_{21}: \cdots x_{2d}:\cdots : x_{(n-1)1}: \cdots : x_{(n-1)d}]\in \mathbb{P}^{d(n-1)-1},
\]
where $p=(z_1,\dots,z_n)\in(\bb{A}^d)^n$, with all $z_i$ distinct, and $z_i-z_n=(x_{i1},\dots,x_{id})\in\bb{A}^d$.

Every compactification of $\tdno$ constructed by Gallardo and Routis is an iterated blow-up of this $\bb{P}^{d(n-1)-1}$ along a collection of collision loci in $\bb{P}^{d(n-1)-1}$ prescribed by a set of weights.

\begin{definition}
    The \emph{domain of admissible weights} is the collection of $n$-tuples of rational numbers :
    \begin{align*}
    \mathcal{D}^T_{d,n} :=
    \left\{
    (a_1, \ldots, a_n) \in \mathbb{Q}^{n} \; | \;
    0 < a_i \leq 1, \textnormal{ for } i\in[n]  \text{ and } a_{[n]}>1 
    \right\}.
    \end{align*}
    We refer to an element $\cal{A}\in\cal{D}^T_{d,n}$ as \emph{collections of weight data}. If two collections of weight data $\cal{A},\cal{B}\in\cal{D}_{d,n}^T$ are such that $a_i\geq b_i$ for all $i\in[n]$, then we write $\cal{A}\geq\cal{B}$.
\end{definition}

Any fixed collection of weight data $\cal{A}\in\cal{D}_{d,n}^T$ uniquely determines a modular compactification $\tdnw$ of $\tdno$\footnote{It is unknown whether these spaces are fine moduli spaces, however, they carry a flat family with $n$ sections, and whose fibers represent the objects being parametrized. This is due to the fact that the spaces are constructed using the Fulton-MacPherson compactification \cite{fulton1994compactification}, for which the same question remains open.}. In order to describe this compactification let us first identify the loci in $\bb{P}^{d(n-1)-1}$ parametrizing configurations where points collide. Define the following subvarieties of $\bb{P}^{d(n-1)-1}$:

\[
    \delta_{ij}=\{q\in\bb{P}^{d(n-1)-1}\,\vert\,x_{ik}=x_{jk}\text{ for all }k\in[d]\}
    \quad\text{and}\quad
    \delta_{in}=\{q\in\bb{P}^{d(n-1)-1}\,\vert\,x_{ik}=0\text{ for all }k\in[d]\}.
\]

More generally, given any $I\subsetneq[n]$ with $|I|\geq 2$, the locus parametrizing configurations $(z_1,\dots,z_n)\in(\bb{A}^d)^n$ where the points $\{z_i\in\bb{A}^d\,\vert\,i \in I\}$ collide is given by
\begin{align}
    \delta_{I}=\left\{
    \begin{aligned}
    \bigcap_{i\neq j\,\in\,I}&\delta_{ij},\text{ if }n\not\in I;\\
    \bigcap_{i\in I\setminus\{n\}}&\delta_{in},\text{ if }n\in I.
    \end{aligned}
    \right.
\end{align}
  
\begin{definition}\label{tdn: def building set}
    Given $\cal{A}\in\mathcal{D}_{d,n}^T$, define the (geometric) building set
    \[
      G_{\cal{A}}=\left\{\delta_{I}\,\vert\, I\subsetneq[n], |I|\geq 2\text{ and }a_I>1\right\}.
    \] 
\end{definition}

\begin{definition}\label{tdn: def tdna}
    Let $\cal{A}\in\mathcal{D}_{d,n}^T$, and endow $G_\cal{A}$ with a total order compatible with containment of subvarieties. Then, the space $\tdnw$ is the wonderful compactification $(\bb{P}^{d(n-1)-1})_{\cal{G}_\cal{A}}$; see Definition~\ref{bg: def wonderful}. In other words, $\tdnw$ is the blow-up of $\bb{P}^{d(n-1)-1}$ along $G_\cal{A}$ in any total order compatible with inclusion of varieties, by Proposition~\ref{bg: li thm 1.3}.
\end{definition}

The set $\cal{D}_{d,n}^T$ admits a coarse wall and chamber decomposition such that $T_{d,n}^\cal{A}$ remains constant within its open chambers. For any $I\subsetneq [n]$, $|I|\geq 2$, the wall $w_I$ is defined by the equation $a_I=1$, and the open chambers are the connected components of $\cal{D}_{d,n}^T\setminus \cup w_I$. It follows from Definitions~\ref{tdn: def building set} and \ref{tdn: def tdna} that $\cal{A},\cal{B}\in\cal{D}_{d,n}^T$ lie in the same coarse chamber if and only if $\cal{G}_\cal{A}=\cal{G}_\cal{B}$, in which case $T_{d,n}^\cal{A}=T_{d,n}^\cal{B}$. In fact, this decomposition can also be obtained from the \emph{fine} chamber decomposition of $\cal{D}_{0,n+1}$ via the inclusion $\cal{D}_{d,n}^T\hookrightarrow\cal{D}_{0,n+1}$ mapping $(a_1,\dots,a_n)\mapsto(a_1,\dots,a_n,1)$.

\begin{proposition}\label{tdn: lemma coarse decompositions agree}
    The coarse chamber decompositions of $\cal{D}_{d,n}^T$ can be obtained via the embedding $\cal{D}_{d,n}^T\hookrightarrow\cal{D}_{0,n+1}$ by intersecting the fine walls of $\cal{D}_{0,n+1}$ with $\cal{D}_{d,n}^T$.
\end{proposition}
\begin{proof}
    The only thing to note is that no fine wall $w_I\subset\cal{D}_{0,n+1}$ with $n+1\in I$ intersects $\cal{D}_{d,n}^T\subset\cal{D}_{0,n+1}$, since $\cal{D}_{d,n}^T=w_{\{n+1\}}\subset\cal{D}_{0,n+1}$.
\end{proof}

Given $\cal{A},\cal{B}\in\cal{D}_{d,n}^T$, let us write $\cal{A}\geq_c\cal{B}$ if there exist collections of weight data $\cal{A}',\cal{B}'\in\cal{D}_{d,n}^T$, in the same coarse chambers as $\cal{A}$ and $\cal{B}$, respectively, and such that $\cal{A}'\geq\cal{B}'$. 

\begin{proposition}{\cite[Proposition~5.1]{gallardo2017wonderful}}
    Consider two collections of weight data $\cal{A},\cal{B}\in\cal{D}_{d,n}^T$ such that $\cal{A}\geq_c\cal{B}$. Then, $\cal{G}_\cal{A}\supseteq\cal{G}_\cal{B}$, and there is a reduction morphism
    \[
        \rho_{\cal{B},\cal{A}}:T_{d,n}^\cal{A}\to T_{d,n}^\cal{B},
    \]
     which is an iterated blow-up along the proper transforms of the elements of $\cal{G}_\cal{A}\setminus\cal{G}_\cal{B}$.
\end{proposition}

In light of the previous discussion, we give the following definitions.

\begin{definition}\label{tdn: def pn}
    Define $\cal{A}_{\bb{P}}^T:=(1/(n-1),\dots,1/(n-1),1/(n-1))\in\cal{D}_{d,n}^T$, so that $(\cal{A}_{\bb{P}}^T)^+=\cal{A}_{\bb{P}}$. From hereafter we consider $\bb{P}^{d(n-1)-1}$ as the compactification $T_{d,n}^\cal{A}$ corresponding to any collection of weight data $\cal{A}\in\cal{D}_{d,n}^T$ defined by the inequalities
    \begin{align}\label{tdn: projective space ineqs}
    \begin{split}
        a_{[n]\setminus\{i\}}\leq 1,\text{ for all }i\in[n].
    \end{split}
    \end{align}
    We denote the interior of this chamber by $Ch_c(\cal{A}_{\bb{P}}^T)$. These inequalities imply that $a_I\leq1$ for all $I\subsetneq[n-1]$.
\end{definition}

In the same spirit, let us define the realizations of the higher-dimensional Losev-Manin space we consider here.

\begin{definition}\label{tdn: def LM}
    Define $\cal{A}_{LM}^T:=(1/(n-1),\dots,1/(n-1),1)\in\cal{D}_{d,n}^T$, so that $(\cal{A}_{LM}^T)^+=\cal{A}_{LM}$. From hereafter we consider $T_{d,n}^{LM}$ as the compactification $T_{d,n}^\cal{A}$ corresponding to any collection of weight data $\cal{A}\in\cal{D}_{d,n}^T$ defined by the inequalities
    \begin{align}\label{tdn: LM ineqs}
    \begin{split}
        &a_{n}+a_i>1,\text{ for all }i\in[n],\\
        &a_{[n-1]}\leq 1.
    \end{split}
    \end{align}
    We denote the interior of this chamber by $Ch_c(\cal{A}_{LM}^T)$. In particular, the second inequality implies that $a_I\leq1$ for all $I\subseteq[n-1]$.
\end{definition}

\begin{remark}
    Strictly speaking $\cal{A}_{\bb{P}}^T$ and $\cal{A}_{LM}^T$ lie on the boundaries of $Ch(\cal{A}_{\bb{P}}^T)$ and $Ch(\cal{A}_{LM}^T)$, respectively. However, by Definition~\ref{tdn: def tdna}, every collection of weight data satisfying the inequalities of Definitions~\ref{tdn: def pn} and \ref{tdn: def LM}, respectively, define the same compactifications.
\end{remark}

If $I\subsetneq[n]$ with $|I|\geq 2$, then $\delta_I$ is torus-invariant with respect to the standard toric structure of $\bb{P}^{d(n-1)-1}$ if and only if $n\in I$. A salient example of this is obtained when considering $\cal{A}=\cal{A}_{LM}^T$, in which case $G_\cal{A}$ is comprised exclusively of torus-invariant subvarieties of $\bb{P}^{d(n-1)-1}$. It follows that all blow-ups in the construction of $T_{d,n}^\cal{A}$ are toric, so the moduli space itself is toric. Following \cite{ggg-higher-lm}, we denote this compactification as $\tdnlm$, and call it the \emph{higher-dimensional Losev-Manin space}.

\begin{example}
    Suppose $d=2$, $n=3$. Then, $T_{2,3}^{\cal{A}}$ is the blow-up of $\bb{P}^{3}$ along $G_\cal{A}$. Let the coordinates of $\bb{P}^{3}$ be $[x_{11}:x_{12}:x_{21}:x_{22}]$.
    
    \begin{itemize}
    \item If $\cal{A}=(1,1,1)$, then $T_{2,3}^{\cal{A}}=T_{d,n}$ and $G_\cal{A} = \left\{\delta_{12},\delta_{13},\delta_{23}\right\}$, where
    \[
    \delta_{12} = \{[x:y:x:y]\},\quad\delta_{13}=\{[0:0:x:y]\},\quad \delta_{23}=\{[x:y:0:0]\}.
    \]
    \item If $\cal{A}=(1/2,1/2,1)$, then $T_{2,3}^{\cal{A}}=\tdnlm$ and $G_\cal{A} = \left\{\delta_{12},\delta_{13}\right\}$,
    because the $a_{1}+a_{2}=1\leq 1$, so $\delta_{12}\not\in G_\cal{A}$. It follows that $T_{2,3}^{LM}$ is the blow-up of $\bb{P}^3$ along two disjoint torus-invariant lines. Note that this is the variety from Examples~\ref{bg: example P3} and \ref{bg: example P3 2}.
    \end{itemize}
\end{example}

\subsection{Proof of Theorem~\ref{tdn: main thm tdn}}

Before presenting the proof of the theorem let us recall its statement.

\begin{theorem}\label{tdn: main thm tdn}[Theorem~\ref{intro: main thm tdn}]
    Let $n\geq 3$ and $\cal{A}\in\cal{D}_{d,n}^T$ be a collection of weight data such that $\cal{A}_{LM}^T\geq_c\cal{A}$ and $\cal{A}\geq_c\cal{A}_{\bb{P}}^T$, so that there are reduction maps $T_{0,n}^{LM}\to T_{d,n}^{\cal{A}}\to\bb{P}^{d(n-1)-1}$. If $n=3$ further suppose that $a_1+a_2\leq 1$. Then,
    \[
        T_{d,n}^{\cal{A}} \cong X(\operatorname{Inf}_d(H_{\cal{A}^+})),
    \]
    where, $\operatorname{Inf}_d(H_{\cal{A}^+})$ is the $d$-inflation of the hypergraph $H_{\cal{A}^+}$ and $X(\operatorname{Inf}_d(H_{\cal{A}^+}))$ its corresponding variety; see Definition~\ref{tdn: def inflation}.
\end{theorem}

As an immediate consequence we obtain the following:

\begin{corollary}
    The higher-dimensional Losev-Manin space $\tdnlm$ is the toric variety corresponding to the $d$-inflation of the $(n-2)$-dimensional permutohedron. 
\end{corollary}

\begin{proof}[Proof of Theorem~\ref{tdn: main thm tdn}]
    Let us recall that $H_{\cal{A}^+}$ is the hypergraph with vertex set $[n-1]$ and hyperedges
    \[
        H_{\cal{A}^+} = \{I\,\vert\, I\subseteq [n-1],|I|\geq 1\text{ and } a_I+a_n>1\}.
    \]
    First consider the case in which $\cal{A}\in Ch_c(\cal{A}_{\bb{P}}^T)$, so that $T_{d,n}^\cal{A}=\bb{P}^{d(n-1)-1}$. In this case $G_\cal{A}=\emptyset$, and $a_I\leq 1$ for all $I\subsetneq[n]$, $|I|\geq 2$. Then, since $a_{[n]}>1$, we obtain that $H_\cal{A}^T=\{[n-1]\}$. It follows that $X(\operatorname{Inf}_d(H_{\cal{A}^+}))=\bb{P}^{d(n-1)-1}$.
    
    Let us now suppose that $T_{d,n}^\cal{A}\ncong\bb{P}^{d(n-1)-1}$. By definition, $X(\operatorname{Inf}_d(H_{\cal{A}^+}))$ is the blow-up of $\bb{P}^{d(n-1)-1}$ along the subvarieties in
    \begin{align*}
        \widehat{G}_\cal{A}:=&\left\{ L_{d,I}:=\bigcap_{i\in I}V(x_{i1},\dots,x_{id})\,\vert\, I\in H_{\cal{A}^+}\setminus\{[n-1]\}\right\}\\
        =&\left\{ L_{d,I}\,\vert\,I\subsetneq [n-1],\, |I|\geq 1, \text{ and }a_I+a_{n}>1\right\}
    \end{align*}
    in order of non-decreasing dimension. On the other hand, since $\cal{A}_{LM}^T\geq\cal{A}$ we have that $\cal{G}_{\cal{A}_{LM}^T}\supseteq\cal{G}_\cal{A}$, so that $T_{d,n}^\cal{A}$ is the blow-up of $\bb{P}^{d(n-1)-1}$ along the subvarieties in
    \[
        G_\cal{A} = \left\{ \delta_{I\cup\{n\}}\,\vert\, I\subsetneq[n-1], |I|\geq 1,\text{ and } a_I+a_n>1\right\}.
    \]
    The reason why only sets of the form $I\cup\{n\}$ appear in this description is because, if $n+1\geq 5$, then $\cal{A}^+\in\cal{D}_{0,n+1}$ and $a_{[n-1]}\leq 1$, by the third clause of Lemma~\ref{mon: lemma chambers pn}. If $n=3$, then $a_1+a_2\leq 1$ by hypothesis. It follows that $\widehat{G}_\cal{A}= G_\cal{A}$, because $L_{d,I} = \delta_{I\cup\{n\}}$ for all $I\subsetneq [n-1]$, $|I|\geq 1$. Moreover, since blow-ups can be performed in any order of non-decreasing dimension in both constructions, one concludes that $X(\operatorname{Inf}_d(H_{\cal{A}})) \cong T_{d,n}^\cal{A}$. 
\end{proof}

\bibliographystyle{plain}
\bibliography{./bibliography}
\end{document}